\documentclass[10pt]{amsart}
\usepackage{amsmath,amssymb,amsthm,latexsym,graphics,psfrag,epsfig}
\usepackage[all]{xy}

\newcommand{\fn}{F_n}
\newcommand{\NN}{\mathcal{N}}
\newcommand{\C}{\mathcal{C}}
\newcommand{\OO}{\mathcal{O}}
\newcommand{\A}{\mathcal{A}}

\newcommand{\Z}{\mathbb{Z}}

\newcommand{\Stab}{\operatorname{Stab}}
\newcommand{\into}{\hookrightarrow}
\newcommand{\Out}{{\rm Out}}
\newcommand{\out}{{\rm Out}(F_n)}
\newcommand{\aut}{{\rm Aut}(F_n)}
\newcommand{\Q}{\mathbb{Q}}
\newcommand{\RR}{\mathfrak{R}}
\newcommand{\kn}{K_n}
\newcommand{\mcg}{MCG^\pm}
\newcommand{\ms}[1][]{(\Sigma_{#1},s_{#1})}
\newcommand{\ums}[1][]{[\Sigma_{#1},s_{#1}]}
\newcommand{\mrg}[1][]{(\Gamma_{#1},g_{#1},\OO_{#1})}
\newcommand{\mg}[1][]{(\Gamma_{#1},g_{#1})}
\newcommand{\rg}[1][]{(\Gamma_{#1},\OO_{#1})}
\newcommand{\pmcg}[3][]{P\Gamma_{#2,#3}^{#1}}
\newcommand{\fmcg}[3][]{\Gamma_{#2,#3}^{#1}}
\newcommand{\direc}[1][e]{\stackrel{\to}{#1}}

\newcommand{\tsigma}{\widetilde\Sigma}

\theoremstyle{plain}
\newtheorem{corollary}{Corollary}
\newtheorem*{corollary*}{Corollary}
\newtheorem{lemma}{Lemma}
\newtheorem*{lemma*}{Lemma}
\newtheorem{theorem}{Theorem}
\newtheorem*{theorem*}{Theorem}
\newtheorem{proposition}{Proposition}

\theoremstyle{definition}
\newtheorem{definition}{Definition}

\theoremstyle{remark}
\newtheorem*{remark}{Remark}
\newtheorem*{warning}{Warning}
\newtheorem{example}{Example}

\newcommand{\bowditechepstein}{MR89e:57004}
\newcommand{\brown}{brown:cohomology}
\newcommand{\cv}{MR87f:20048}
\newcommand{\harer}{MR87f:57009}
\newcommand{\harerimproved}{harerimproved}
\newcommand{\harerexact}{harerimproved}
\newcommand{\harervcd}{MR87c:32030}

\newcommand{\hatvogtstab}{math.GT/0406377}
\newcommand{\hoare}{MR80e:20050} 
\newcommand{\horak}{horak}  
\newcommand{\ivanov}{MR91g:57010}
\newcommand{\madsenweiss}{math.AT/0212321}
\newcommand{\mccooltwo}{MR53:624}
\newcommand{\mulasepenkava}{MR2001g:30028}
\newcommand{\mks}{MR54:10423}
\newcommand{\penner}{MR89h:32044} 
\newcommand{\quillen}{MR49:2895}
\newcommand{\vogtsurvey}{MR1950871}

\newcommand{\zies}{MR82h:57002}

\begin{document}

\title{Mapping class subgroups of $\out$}
\author{Matthew Horak}
\maketitle

\begin{abstract}
We construct a covering of the spine of Culler-Vogtmann Outer space by complexes of ribbon graphs. By considering the equivariant homology for the action of $\out$ on this covering, we construct a spectral sequence converging to the homology of $\out$ that has its $E^1$ terms given by the homology of mapping class groups and their subgroups.  This spectral sequence can be seen as encoding all of the information of how the homology of $\out$ is related to the homology of mapping class groups and their subgroups
\end{abstract}

\section{Introduction}\label{S:intro}

Much is known about the cohomology of mapping class groups of surfaces.  Throughout this paper, we will consider only orientable surfaces.  Thus, even when not explicitly stated, all surfaces will be orientable.  Let $\Sigma$ be an orientable surface with boundary, and let $P\Gamma(\Sigma)$ be the group of isotopy classes, relative to the boundary, of homeomorphisms of $\Sigma$ that fix the boundary.  We call $P\Gamma(\Sigma)$ the pure mapping class group of $\Sigma$.  In~\cite{\harer}, Harer proved that the  $k^{th}$ integral (and therefore also rational) homology group of $P\Gamma(\Sigma)$ is independent of the genus and number of boundary components of $\Sigma$ if the genus of $\Sigma$ is at least $3k$. Later, Ivanov~\cite{\ivanov} and Harer~\cite{\harerimproved} improved these bounds, and Harer was able to find the exact location at which the rational homology stabilizes~\cite{\harerexact}.   Harer~\cite{\harervcd} has also computed the virtual cohomological dimension of $P\Gamma(\Sigma)$ and shown that this group has no rational homology at its VCD.  Madsen and Weiss~\cite{\madsenweiss} have determined the entire stable integral cohomology algebra of pure mapping class groups.  In particular, their result verifies the conjecture of Mumford that the stable rational cohomology algebra is a polynomial algebra with a single generator in each even dimension.

For outer automorphism groups of free groups, much less is known.  In ~\cite{\cv}, Culler and Vogtmann compute the VCD of $\out$ by considering the action of this group on a contractible simplicial complex known as the spine of outer space.  Recently, Hatcher and Vogtmann~\cite{\hatvogtstab} have showen that the $k^{th}$ integral homology of $\out$ is independent of $n$ if $n \geq 2k+5$, but the exact stability range remains unknown.  Indeed, there are no nontrivial stable rational homology or cohomology classes known for $\out$.  For a good survey of the current state of knowledge about $\out$ and $\aut$, see~\cite{\vogtsurvey}.

There are many varieties of mapping class groups, all of which are related to each other and to $\out$.  As mentioned above, for a surface $\Sigma$ with boundary, the pure mapping class group of $\Sigma$ consists of isotopy classes, relative to the boundary, of diffeomorphims that fix the boundary pointwise.  Thus, pure mapping classes of a bounded surface preserve orientation.  For punctured surfaces, the mapping class group is simply the group of isotopy classes of orientation preserving diffeomorphisms of the surface.  The group of all isotopy classes of diffeomorphisms of a punctured surface will be called the extended mapping class group of that surface, so that extended mapping class groups contain orientation reversing diffeomorphisms while ordinary mapping class groups do not.

These classes of mapping class groups are related by homomorphisms induced by inclusion maps of the associated surfaces.  Given a bounded surface, we can produce a punctured surface with the same genus and fundamental group by gluing a punctured disk to each boundary component.  The inclusion of a bounded surface into a punctured surface given by this gluing induces a homomorphism from the mapping class of the bounded surface into the mapping class group of the punctured surface.  The kernel of this map is generated by the Dehn twists about curves parallel to the boundary.  As for mapping class groups and extended mapping class groups of punctured surfaces without boundary, the mapping class group is simply a index two subgroup of the extended mapping class group.

The relation of mapping class groups of surfaces to $\out$ comes from a theorem of Zieschang~\cite[Theorem 5.15.3]{\zies}, which says that for a punctured surface, $\Sigma$ the mapping class group of $\Sigma$ is isomorphic to the subgroup of $\Out(\pi_1(\Sigma))$ that stabilizes the set conjugacy classes in $\pi_1(\Sigma)$ represented by simple closed curves around the punctures.  For extended mapping class groups, we must take the subgroup of outer automorphisms that either permutes the above conjugacy classes or permutes them and inverts them all, which is the same as the stabilizer of the set of these conjugacy classes and their inverses.  Since $\pi_1(\Sigma)$ is free, this theorem gives the mapping and extended mapping class groups of punctured surfaces as subgroups of $\out$.

This paper represents an attempt to clarify the relationship between the homology $\out$ and that of mapping class groups of surfaces.  We construct a first quadrant spectral sequence converging to $H_*(\out)$.  The spectral sequence arises from a covering of the spine of Outer space by a collection of subcomplexes called ribbon graph subcomplexes.  We prove that the nerve of this covering is contractible.  The spectral sequence mentioned above is the equivariant homology spectral sequence of the action of $\out$ on this nerve.

All of the terms on the $E^1$ page of this spectral sequence are given by the homology simplex stabilizers.  For a $0$-simplices, the stabilizer is simply extended mapping class group of a punctured surface $\Sigma$, or equivalently the stabilizer of the set conjugacy classes in $\fn$ that correspond to positively and negatively oriented curves about the punctures of $\Sigma$.  For higher dimensional simplices, stabilizers are given by the generalized stabilizers, $\A_{U,G}$, of $n$-tuples of conjugacy classes, which are studied by McCool~\cite{\mccooltwo}.  Equivalently, these groups are finite index subgroups of the stabilizers of certain sets of conjugacy classes in $\fn$.  We prove,
\begin{theorem*} For any $\out$-module, $M$, there is a spectral sequence of the form
\begin{equation*}
E^1_{pq} = \bigoplus_{\sigma\in\Delta_p}H_q(G_\sigma;M_\sigma)\Rightarrow H_{p+q}(\out;M),
\end{equation*}
where $\Delta_0$ is the set of homeomorphism classes of punctured orientable surfaces with fundamental group $\fn$ and for vertex $v\in\Delta_0$ corresponding to surface $\Sigma$, the stabilizer $G_v$ is the extended mapping class group $\mcg(\Sigma)$.  Moreover, for $p>0$, each $G_\sigma$ is a generalized stabilizer of the form $\A_{U_\sigma,H_\sigma}$.
\end{theorem*}

The rest of this paper is organized as follows.  In Sections $2$ and $3$ we review the definitions of Outer space, the spine of Outer space, ribbon graphs and some related objects.  In Section $4$, we construct a covering of the spine of Outer space by subcomplexes of ribbon graphs.  Section $5$ is devoted to the proof of the fact that the nerve of this covering is contractible.  In Section $6$ we determine simplex stabilizers for the action of $\out$ on the nerve.  The analysis of the equivariant homology spectral sequence for this action appears in final two sections where we prove the above theorem and use Harer's stability theorems to find rough upper bounds on the dimensions of some portions of the $E^\infty$ page of the spectral sequence.  Harer's stability theorems are stated for surfaces with boundary but, until we need them in Section 8, we will focus strictly on surfaces with punctures.  In that section we will adapt Harer's results to our situation.

\section{Outer space}
For convenience and to set notation, we briefly review the construction in~\cite{\cv} of Outer space and its spine.  A {\em graph} is a connected, 1-dimensional CW-complex.  We will consider only finite graphs with all vertices having valence at least $3$.  A {\em subforest} of a graph $\Gamma$ is a subgraph of $\Gamma$ that contains no circuits; a forest is a disjoint union of trees.

Fix an integer $n \geq 2$.  Denote by $R_0$ the standard $n$-petal rose; $R_0$ has one vertex and $n$ edges.  Fix an identification $\pi_1(R_0) = \fn$.  A {\em marking} on a graph is a homotopy equivalence,  $g:R_0 \rightarrow \Gamma$.  We define an equivalence relation on the set of markings by $(\Gamma_1,g_1) \sim (\Gamma_2,g_2)$ if there is a graph isomorphism $h:\Gamma_1\rightarrow\Gamma_2$ making the following diagram commute up to free homotopy:
\begin{equation}\label{E:eq.relation}
\xymatrix@C-20pt@R-20pt{
& & \Gamma_1\ar[dd]^h\\
R_0\ar[rru]^{g_1} \ar[rrd]_{g_2} & &\\
& & \Gamma_2.
}
\end{equation}
That is, $g_2 \simeq h \circ g_1$.  An equivalence class of markings is called a {\em marked graph} and also denoted $\mg$.  The marking $g$ identifies $\pi_1(\Gamma)$ with $\fn$ up to composition with an inner automorphism.

If $\Phi$ is a forest in the marked graph $\mg$, then collapsing each component of $\Phi$ to a point produces another marked graph, denoted $(\Gamma / \Phi , q \circ g)$, where $q$ is the quotient map collapsing each component of $\Phi$ to a point.  Passing from $\mg$ to $(\Gamma / \Phi , q \circ g)$ is called a {\em forest collapse}.  There is a partial order on the set of marked graphs with fundamental group $\fn$ defined by $\mg[1] \leq \mg[2]$ if there is a forest collapse taking $\mg[2]$ to $\mg[1]$.  The geometric realization of the poset of marked graphs is the spine of Outer space and is denoted by $K_n$.

The group $\out$ acts on $K_n$ by changing the markings of the underlying graphs.  Explicitly, for $\psi \in \out$, 
\begin{equation}\label{eq:out.action}
    \mg \cdot \psi := (\Gamma, g\circ |\psi|),
\end{equation}
where $|\psi|: R_0 \to R_0$ is a homotopy equivalence inducing an automorphism of $\fn = \pi_1(R_0)$ that represents the outer automorphism class $\psi$.  Culler and Vogtmann observe that this action is cocompact and that vertex stabilizers are finite.

Culler and Vogtmann also define a larger space, Outer space, consisting of metric marked graphs.  This space has the disadvantage of not being a simplicial complex and the $\out$ action not being cocompact.  The complex $\kn$ can be constructed as a simplicial spine onto which of Outer space deformation retracts.

\section{Ribbon graphs}
\begin{figure}
    \psfrag{t}{\Huge $\rightarrow$}
    \psfrag{h}{\huge $\approx$}
    \psfrag{g1}{\large $(\Gamma,\OO_1)$}
    \psfrag{g2}{\large $(\Gamma,\OO_2)$}
    \resizebox{5in}{2in}{\includegraphics{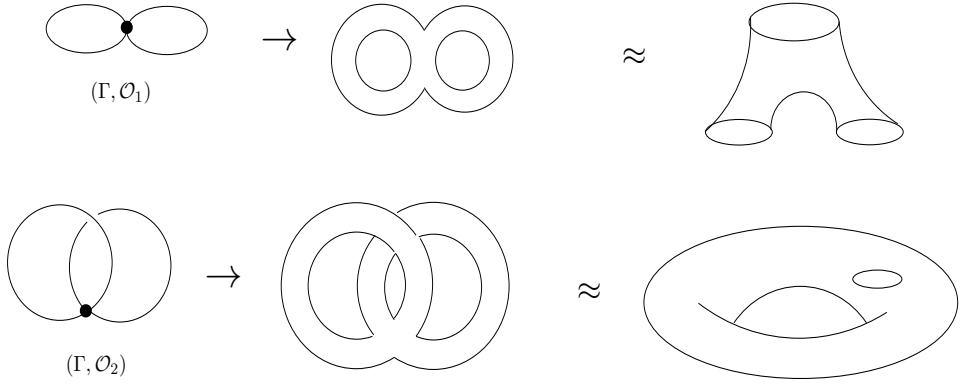}}
    \caption{Fattenings of ribbon graphs}\label{F:fatten}
\end{figure}

There are similar constructions for mapping class groups that use marked ribbon graphs rather than ordinary marked graphs.  A {\em ribbon graph} is a graph $\Gamma$ together with, at each vertex $v$, a cyclic ordering of the set $h(v)$ of half edges incident to $v$.  The collection of cyclic orderings at the vertices is called a {\em ribbon structure} for $\Gamma$, and is denoted by $\OO$.  The term ``ribbon graph'' is used because one can construct a bounded surface from a ribbon graph $\rg$ by fattening its edges to ribbons.  We give a formal construction of this surface after Definition~\ref{D:bd.cycle}, but informally, the surface is constructed from $\rg$ by replacing each edge by a ribbon and gluing the ribbons together at their ends according to the cyclic order of the corresponding half edges.  The gluing is done in such a way as to produce an oriented surface.  Figure~\ref{F:fatten} shows this process for two different ribbon structures on a rose with $2$ edges.  In these figures, ribbon structures are specified by the given embeddings of a neighborhood of the vertices into the plane.  The ribbon graph $(\Gamma,\OO_1)$ produces a pair of pants while $(\Gamma,\OO_2)$ produces a torus with one boundary component.

The boundary curves of the surface produced from $\rg$ correspond to reduced edge paths in $\Gamma$ that follow the cyclic ordering at the vertices in the sense of the following definition.  Following~\cite{\mulasepenkava}, we view a directed edge $\direc$ as an ordering $(e^+,e^-)$ of the half edges $e^+$ and $e^-$ comprising $e$.  
\begin{definition}\label{D:bd.cycle}A {\em boundary cycle} in the ribbon graph $\rg$ is a directed reduced edge cycle,
$$(\direc[e_1],\direc[e_2],\ldots,\direc[e_{l-1}],\direc[e_l]=\direc[e_1])$$
such that for each $i$ the half-edges $e_i^+$ and $e_{i+1}^-$ are incident to the same vertex, and in the cyclic ordering at that vertex, $e_{i+1}^-$ directly follows $e_i^+.$
\end{definition}

For our purposes, it will be more convenient to work with punctured surfaces, so we now give a precise construction of a punctured surface $|\Gamma,\OO|$ from a ribbon graph $\rg$.  First note that each edge of $\Gamma$ is traversed exactly once in each direction by the set of boundary cycles of $\rg$.  Construct a space $|\Gamma,\OO|$ by gluing a once punctured disk to $\Gamma$ along each boundary cycle $\gamma$ of $\rg$.  By verifying that a small neighborhood of each vertex in $|\Gamma,\OO|$ is indeed a disk, one can verify that $|\Gamma,\OO|$ is a surface that deformation retracts onto $\Gamma$.  One can also verify that $|\Gamma,\OO|$ is orientable and we orient it such that a small positively oriented simple closed curve around a vertex $v$ of $\Gamma$ intersects the half edges in $h(v)$ in the cyclic order determined by $\OO$.

If $\Gamma$ is marked by the homotopy equivalence $g: R_0 \to \Gamma$, then the composition of $g$ with the inclusion $i: \Gamma \into |\rg|$ is a homotopy equivalence that identifies $\pi_1(\Sigma)$ with $\fn$ up to inner automorphism, just as in the case of marked graphs.  This gives the notion of a homotopy marked surface.

\begin{definition}
A {\em homotopy marked surface} is an equivalence class of pairs $\ms$, where $\Sigma$ is a punctured, orientable surface with $\pi_1(\Sigma) \cong \fn$ and $s: R_0 \to \Sigma$ is a homotopy equivalence.  The equivalence relation on pairs is given by $\ms[1] \sim \ms[2]$ if there is an orientation preserving homeomorphism $h: \Sigma_1 \to \Sigma_2$ with $h\circ s_1 \simeq s_2$.
\end{definition}

Recall that we have fixed an integer $n\geq 2$.  Often we drop the word ``homotopy" and simply use ``marked surface" for a homotopy marked surface.  Unless otherwise stated, marked surfaces will always be of type $\Sigma_g^s$ with $ s>0$ and $2g+s-1=n$;  that is they will be punctured surfaces without boundary and with fundamental group $\fn$.  In Section~\ref{s:cover}, the equivalence relation defined by homeomorphisms that do not necessarily preserve orientation will be useful.  We will denote this equivalence relation by $\sim_\pm$, and use square brackets to denote its equivalence classes, $\ums$.  We say that the marked graph $\mg$ can be {\em drawn in} the marked surface $\ms$ is there is a ribbon structure $\OO$ on $\Gamma$ such that $|\mrg| \sim \ms$.  In this case, there is an embedding $i: \Gamma \hookrightarrow \Sigma$ such that $s \simeq i \circ g$.

\begin{definition}
The {\em ribbon graph complex} for the marked surface $\ms$ is the subcomplex of $\kn$ spanned by graphs that can be drawn in $\ms$.  This complex is denoted by $\RR_{\ms}$.
\end{definition}

We will often identify a marked graph or ribbon graph with the corresponding vertex of $\kn$ or $\RR_{\ms}$.  Thus, for example, if $\rho$ is a marked rose in $\RR_{\ums}$ then $lk_{\RR_{\ms}}(\rho)$ will be the link in $\RR_{\ms}$ of the vertex corresponding to $\rho$.

The ribbon graph complex $\RR_{\ms}$ and related complexes have been important tools in the study of mapping class groups surfaces.  In particular, $\RR_{\ms}$ is a subcomplex of the first barycentric subdivision of the arc complex that Harer uses compute the VCD of the pure mapping class group of a surface with boundary~\cite{\harervcd}.  Also, for a punctured surface $\Sigma$, Bowditch and Epstein~\cite{\bowditechepstein} and Penner~\cite{\penner} use arc systems on $\Sigma$ to give an open cell decomposition of a space they call the decorated Teichm\"uller space of $\Sigma$.  By taking the dual graph of an arc system in $\Sigma$, this decomposition may be interpreted in terms of metric ribbon graphs.  In the same way that $\kn$ is a simplicial spine of Outer space, $\RR_\Sigma$ may be seen as a simplicial spine of the decorated Teichm\"uller space of $\Sigma$.

\section{The ribbon cover of $\kn$}\label{s:cover}

We first note that the ribbon subcomplex of $\kn$ associated to a marked surface does not depend on the orientation taken on the surface.  This is because if the ribbon structure $\OO$ draws $\mg$ in $\ms$, then $\OO^{op}$ draws $\mg$ in $\ms^{op}$, where $\OO^{op}$ is the ribbon structure obtained by reversing all cyclic ordering of $\OO$ and $\ms^{op}$ the marked surface obtained by reversing the orientation of $\ms$.  Therefore, there is a well-defined subcomplex  $\RR_{\ums}$ of $\kn$.  We have,

\begin{lemma}\label{L:kn.covered}
$\kn$ is covered by its ribbon graph subcomplexes.
\end{lemma}
\begin{proof}
Recall that $\kn$ is the geometric realization of the poset of marked graphs with fundamental group $\fn$, so the vertices of $\kn$ are partially ordered.  For a vertex $v$ of $\kn$, let $st(v)$ be the star of $v$ and let $st^+(v)$ denote the subcomplex of $st(v)$ spanned by $v$ together with vertices of $st(v)$ that are greater than $v$ in the partial order.  Similarly let $st^-(v)$ denote the subcomplex of $st(v)$ spanned by $v$ and vertices less than $v$.  Thus, if $v$ corresponds to the marked graph $\mg$, then $st^+(v)$ consists of the vertices of $\kn$ corresponding to graphs that may be collapsed to $\mg$ and $st^-(v)$ consists of vertices corresponding to graphs to which $\mg$ collapses.

Now, suppose that $v \in \RR_{\ums}$ corresponds to the marked graph $\mg$.  Then $\Gamma$ has a ribbon structure $\OO$ that draws $\mg$ in $\ms$.  If $e$ is any edge in $\Gamma$ that is not a loop, then the marked graph $(\Gamma / e, q \circ g)$ inherits a ribbon structure $\OO / e$ from $\OO$ that draws $(\Gamma / e , q \circ g)$ in $\ms$.  Therefore, $st^-(v) \subset \RR_{\ums}$.

To see that every simplex of $\kn$ belongs to some ribbon graph subcomplex, let $\sigma$ be a simplex of $\kn$.  If $w$ is the vertex of $\sigma$ that is the greatest in the partial ordering of the vertices, then $\sigma$ is contained in the complex $st^-(w)$.  Suppose that $w$ corresponds to the marked graph $\mg[0]$  Choose any ribbon structure $\OO_0$ on $\Gamma_0$ and set $\ms := |\mrg[0]|$.  Then $w \in \RR_{\ums}$ so that $st^-(w)$ is contained in $\RR_{\ums}$.  Since $\sigma$ has $w$ as its greatest vertex, $\sigma$ is a simplex of $st^-(w)$, which is contained in $\RR_{\ums}$.  Therefore every simplex of $\kn$ is contained in some ribbon graph complex.
\end{proof}

We begin our study of the nerve of this cover with several definitions and lemmas.

\begin{definition}\label{D:bd.classes}
Suppose that the homotopy marked, oriented surface $\ms$ has $k$ punctures, $p_1,\ldots,p_k$.  Let $\gamma_j$ be a simple closed curve in $\Sigma$ that disconnects $\Sigma$ by cutting off a disk punctured at $p_j$.  By virtue of the marking and orientation of $\Sigma$, $\gamma_j$ corresponds to a conjugacy class in $F_n$.  The set of such conjugacy classes is called the set of {\em boundary classes} of $\Sigma$ and denoted $W_{\ms}$ or simply $W_\Sigma$.  Similarly, the set of conjugacy classes in $\fn$ represented by the boundary cycles of the marked ribbon graph $\mrg$ is called the set of {\em boundary classes} of $\mrg.$  Note that if $\ms =|\mrg|$ then the boundary classes of $\ms$ and the boundary classes of $\mrg$ are the same.
\end{definition}

\begin{lemma}\label{L:unique.ribbon}  If the marked graph $\mg$ can be drawn in $\ms$, then $\mg$ has exactly one ribbon structure giving $\ms$.
\end{lemma}

\begin{proof}
Since $\Gamma$ can be drawn in $\ms$, there is a ribbon structure $\OO$ on $\mg$ with $|\mrg| = \ms$.  Suppose that $\OO'$ is a different ribbon structure on $\Gamma$.  Since $\OO'$ is different from $\OO$, we may choose a vertex $v$ and half edges $e^+, e_1^-$ and $e_2^-$ of $\Gamma$ with $e_1^- \neq e_2^-$, with $e_1^-$ following $e^+$ in the cyclic ordering $\OO$ but with $e_2^-$ following $e^+$ in the cyclic ordering $\OO'$.  This means that the sequence $\direc[e]\direc[e_1]$ appears in the boundary cycles of $(\Gamma,\OO)$ while the sequence $\direc[e]\direc[e_2]$ appears in the boundary classes of $(\Gamma,\OO')$.  Since each directed edge of $\Gamma$ appears exactly once in the set of boundary cycles for any given ribbon structure, the set of boundary cycles of $(\Gamma,g,\OO')$ must differ from those of $\mrg$.  Therefore, the set of boundary classes of $(\Gamma,g,\OO')$ differ from those of $\mrg$.  By the above remark $|(\Gamma,g,\OO')|$ must be different from $\ms$ because equivalent marked surfaces have the same boundary classes.
\end{proof}

\begin{remark} Each of the two orientations of $\Sigma$ gives a unique ribbon structure to $\mg$.  These two ribbon structures are opposite of each other.
\end{remark}

We now prove that $\ms$ and $\ms[1]$ give the same ribbon graph subcomplexes of $\kn$ if and only if $\ums = \ums[1]$.  The main step is the following lemma.

\begin{lemma}\label{L:construct.order}  Let $\ms$ be a marked surface and let $\RR = \RR_{\ums}$ be the corresponding ribbon graph subcomplex of $K_n$. The ribbon structure given by $\ms$ to a marked rose $\rho$ in $\RR$ can be reconstructed, up to reversal of the cyclic order, by the (non-ribbon) marked graphs in $lk_{\mathfrak R}(\rho)$.
\end{lemma}

\begin{proof}
Choose a direction for each edge in $\rho$.  The marking of $\rho$ determines a labelling of the directed edges by a basis $X = \{a_1,a_2,\ldots,a_n\}$ of $\fn$.  The basis $X$ and this labelling are determined up to composition with an inner automorphism of $\fn$.

Consider the marked graphs in $lk_{\kn}(\rho)$ with exactly two vertices, one of which is trivalent.  These graphs are constructed from $\rho$ as follows.  Let $e^+$ and $f^+$ be any two half edges of $\rho$.  Construct a new marked graph $\rho(e^+,f^+)$ by deleting the vertex of $\rho$ and replacing it with two new vertices $v_0$ and $v_1$ joined by a new edge $\tilde{e}$.  Attach the half edges $e^+$ and $f^+$ to $v_0$ and attach the rest of the half edges of $\rho$ to $v_1$.  Mark $\rho(e^+,f^+)$ so that collapsing $\tilde{e}$ to a point gives the original marking on $\rho$.

By Lemma~\ref{L:unique.ribbon}, $\rho$ has exactly one ribbon structure giving $\ms$.  This means that if we allow for orientation reversing homeomorphisms of the surface, $\rho$ has two ribbon structures giving $\ums$, and these ribbon structures are opposite of each other.  A marked graph of the form $\rho(e^+,f^+)$ lies in $lk_{\RR}(\rho)$ if and only if $e^+$ and $f^+$ are adjacent in these ribbon structures.  Thus, given a half edge, $a_i^+$ in $\rho$, exactly two graphs of the form $\rho(a_i^+,a_j^{\epsilon_j})$ and $\rho(a_i^+,a_k^{\epsilon_k})$ will lie in $lk_{\RR}(\rho)$ and they will be the graphs for which $a_i^+$ is adjacent to the half edges $a_j^{\epsilon_j}$ and $a_k^{\epsilon_k}$ in the ribbon structure on $\rho$.  Therefore, for each $i$, the non-ribbon graphs in $lk_{\RR}(\rho)$ determine the half edges adjacent to $a_i^+$ and $a_i^-$ in the ribbon structure on $\rho$.  There are only two cyclic orderings of the half-edges that satisfy this adjacency data, and they are opposites of each other.  Example~\ref{EX:cyc.order} works this out for a ribbon rose with $n=3$.
\end{proof}

\begin{example}\label{EX:cyc.order}
Figure~\ref{F:cyc.order}, shows some of the graphs in the link of a marked ribbon rose in $\RR_\Sigma$. The fact that graphs $(1)$ and $(6)$ have the half-edges $a^+$ and $c^+$, respectively, adjacent to the edge $a^-$ implies that the half-edges in $\rho$ adjacent to $a^-$ are $a^+$ and $c^+$.  Similarly, the half-edges adjacent to any other half-edge can be determined by some pair of the graphs in Figure~\ref{F:cyc.order}.

\begin{figure}

\begin{center}
\begin{minipage}{.8 in}

\psfrag{a}{\Huge $a$}
\psfrag{b}{\Huge $b$}
\psfrag{c}{\Huge $c$}
\psfrag{r}{\Huge $\rho$}
\resizebox{.7 in}{.7 in}{\includegraphics{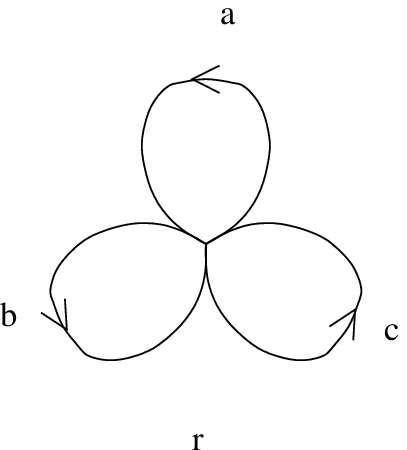}}\\

\end{minipage}
\end{center}

\begin{center}
\begin{minipage}{3.08 in}
\begin{minipage}{1 in}
\begin{center}
\psfrag{a}{\Huge $a$}
\psfrag{b}{\Huge $b$}
\psfrag{c}{\Huge $c$}
\resizebox{.5 in}{.7 in}{\includegraphics{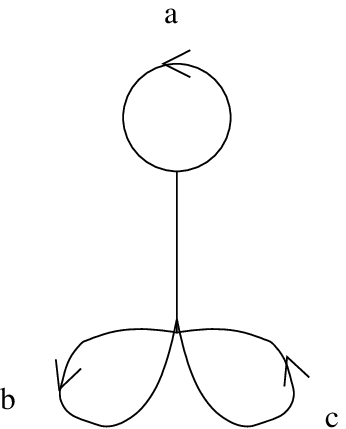}}\\
$(1)$

\vspace{.2 in}

\psfrag{a}{\huge $a$}
\psfrag{b}{\huge $b$}
\psfrag{c}{\huge $c$}
\resizebox{.5 in}{.7 in}{\includegraphics{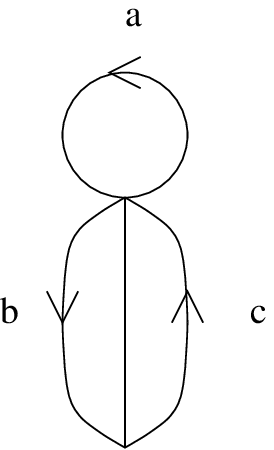}}\\
$(4)$
\end{center}
\end{minipage}
\begin{minipage}{1 in}
\begin{center}
\psfrag{a}{\huge $c$}
\psfrag{b}{\huge $a$}
\psfrag{c}{\huge $b$}
\resizebox{.5 in}{.7 in}{\includegraphics{blowrose2.eps}}\\
$(2)$

\vspace{.2 in}

\psfrag{a}{\Huge $c$}
\psfrag{b}{\Huge $a$}
\psfrag{c}{\Huge $b$}
\resizebox{.5 in}{.7 in}{\includegraphics{blowrose1.eps}}\\
$(5)$
\end{center}
\end{minipage}
\begin{minipage}{1 in}
\begin{center}
\psfrag{a}{\Huge $b$}
\psfrag{b}{\Huge $c$}
\psfrag{c}{\Huge $a$}
\resizebox{.5 in}{.7 in}{\includegraphics{blowrose1.eps}}\\
$(3)$

\vspace{.2 in}

\psfrag{a}{\huge $b$}
\psfrag{b}{\huge $c$}
\psfrag{c}{\huge $a$}
\resizebox{.5 in}{.7 in}{\includegraphics{blowrose2.eps}}\\
$(6)$
\end{center}
\end{minipage}
\end{minipage}
\end{center}

    \caption{Some Graphs in $lk^+_{\RR_\Sigma}(\rho)$}\label{F:cyc.order}
\end{figure}
\end{example}

\begin{proposition}\label{P:equivalence}
$\RR_{\ums[1]}=\RR_{\ums[2]}$ if and only if $\ums[1] = \ums[2]$.
\end{proposition}

\begin{proof}
First suppose that $\ms[1]$ and $\ms[2]$ are equivalent via the (possibly orientation-reversing) homeomorphism $h:\Sigma_1\rightarrow\Sigma_2$. Then $h$ can be used to draw in $(\Sigma_2, s_2)$ any graph that can be drawn in $(\Sigma_1, s_1)$ and $h^{-1}$ can be used to draw in $(\Sigma_1, s_1)$ any graph that can be drawn in $(\Sigma_2, s_2)$, so $\RR_{\ums[1]} = \RR_{\ums[2]}$.

Now, suppose that $\RR_{\ums[1]}=\RR_{\ums[2]}$.  Let $\RR := \RR_{\ums[1]}=\RR_{\ums[2]}$.  Fix a marked rose, $\rho\in\RR$.  The rose $\rho$ inherits a ribbon structure from $\RR_{\ums[1]}$ that gives $\ms[1]$ and a ribbon structure from $\RR_{\ums[2]}$ that gives $\ms[2]$.  By Lemma~\ref{L:construct.order}, these structures are determined up to reversal by the non-ribbon graphs in $lk_{\RR_{\ums[1]}}(\rho)$ and $lk_{\RR_{\ums[2]}}(\rho)$ respectively.  But, $\RR_{\ums[1]}=\RR_{\ums[2]}$ so $lk_{\RR_{\ums[1]}}(\rho) = lk_{\RR_{\ums[2]}}(\rho)$.  Therefore, the ribbon structures must coincide or be opposites of each other.  In the first case, $(\Sigma_1,s_1) \sim (\Sigma_2,s_2)$ and in the second case, $(\Sigma_1,s_1) \sim (\Sigma_2,s_2)^{op}$.  Thus, $\ums[1] = \ums[2]$.
\end{proof}

This proposition gives a convenient description of the covering of $\kn$ by its ribbon graph subcomplexes.  The covering is locally finite because each different homotopy marked surface that contains a specific graph endows that graph with a different ribbon structure.  A graph has only finitely many different ribbon structures so a given marked graph can be drawn in only finitely many marked surfaces and hence lies in only finitely many different ribbon graph subcomplexes.

Let $\NN_n$ denote the nerve of this covering.  That is, $\NN_n$ is the simplicial complex containing a $k$-simplex, $\langle \RR_{\ums[0]},\cdots,\RR_{\ums[k]} \rangle$ for every collection of ribbon graph complexes, $\{\RR_{\ums[0] },\cdots,\RR_{\ums[k]}\}$, such that the intersection $\bigcap_{i=0}^{k} \RR_{\ums[i]}$ is nonempty.  By Proposition~\ref{P:equivalence}, the vertex set of $\NN_n$ is the set of unoriented equivalence classes, $\ums$, of marked surfaces.

The action of $\out$ on $\kn$ permutes the ribbon graph subcomplexes because if $\mg$ can be drawn in $\ms$, then $\mg \cdot \psi = (\Gamma, g \circ |\psi|)$ can be drawn in $(\Sigma, s \circ |\psi|)$.  Therefore, $\out$ maps intersections of ribbon graph subcomplexes to intersections of ribbon graph subcomplexes, so it acts on $\NN_n$.  The equivariant homology of this action provides a spectral sequence that relates the homology of $\out$ to that of mapping class groups.

Although it will not be necessary for the development here, we remark briefly on the compactness properties of $\NN_n$ and the $\out$ action.  In general ($n\geq 3$), all vertices of $\NN_n$ have infinite valence because, for $n\geq 3$, the ribbon complex for any homotopy marked surface intersects the ribbon complexes of infinitely many other homotopy marked surfaces, though not all in the same place.  This is because the ribbon graph subcomplex of any marked surface $\ums$ with fundamental group of rank at least $3$, contains infinitely many different marked roses, which can be seen as follows.  If $\rho$ is a marked rose with edges labelled by the basis, $X = \{a_1, \ldots, a_n\}$, then $\rho$ can be drawn in a marked $(n+1)$-times punctures sphere $\Sigma_1$ with boundary classes, $$W_{\Sigma_1} = \{a_1, \ldots, a_n, a_n^{-1}\cdots a_1^{-1}\}.$$  Since two marked spheres with different boundary classes cannot be equivalent, the infinitely many different marked roses in $\RR_{\ums}$ give rise to infinitely many different marked spheres all of whose ribbon complexes intersect $\RR_{\ums}$.  Therefore, the vertex of $\NN_n$ corresponding to $\ums$ has infinite valence.  On the other hand we have,

\begin{proposition}\label{P:cocompact}
$\out$ acts cocompactly on $\NN_n$.
\end{proposition}

\begin{proof}
Fix a marked rose $\rho$.  For each $p$-simplex $\langle \Sigma_0,\ldots,\Sigma_p \rangle$, the subcomplex $\bigcap_{i=0}^{p} \RR_{\Sigma_i}$ contains a rose.  This rose may be taken to $\rho$ by an element of $\out$, so each orbit of $p$-simplex has a representative all of whose surfaces contain $\rho$.  Since a marked rose can be drawn in only finitely many different marked surfaces, there are only finitely many orbits of $p$-simplices.
\end{proof}

\section{Contractibility of $\NN_n$}
To show that $\NN_n$ is contractible, we will need the following result from $\check{{\rm C}}$ech theory.

\begin{lemma}\label{L:nerve.contract}
Let $\mathfrak U$ be a cover of the CW-complex $X$ by a family of subcomplexes.  If every nonempty intersection of finitely many complexes in $\mathfrak U$ is contractible, then the nerve of the cover is homotopy equivalent to $X$.
\end{lemma}

We will apply Lemma~\ref{L:nerve.contract} to the covering of $\kn$ by ribbon graph complexes.  Thus, the rest of this section is devoted to the proof of,

\begin{proposition}\label{P:intersections} For any finite collection $\{\RR_{\Sigma_0},\ldots,\RR_{\Sigma_k}\}$ of ribbon graph subcomplexes of $K_n$, the subcomplex
$$\bigcap_{i=0}^k \RR_{\Sigma_i}$$
of $\kn$ is either empty or contractible.
\end{proposition}

The main tool in analyzing these intersections will be the $K_{min}$ subcomplexes of $K_n$, which are used by Culler and Vogtmann~\cite{\cv} to show that $\kn$ is contractible.  The definition of the $K_{min}$ complexes involves the following norm defined for each set of conjugacy classes in $\fn$.  Let $\C$ denote the set of all conjugacy classes of elements in $\fn$.  For a marked rose $\rho$, an element $w\in\C$ can be represented by a unique reduced edge path in $\rho$.

\begin{definition}\label{D:norm}  Let $W$ be a finite set of conjugacy classes of $\fn$ and $\rho$ a marked rose in $\kn$.  The {\em norm}, $\|\rho\|_W$, of $\rho$ with respect to $W$ is the sum of the number of edges in each reduced edge path in $\rho$ that corresponds to an element of $W$.
\end{definition}

If $X$ is a basis labelling the edges of $\rho$, then $\|\rho\|_W$ is sometimes written $\|X\|_W.$  The $K_{min}$ subcomplex for $W$ is defined as the union of the stars of the roses $\rho$ for which $\|\rho\|_W$ is minimal over all marked roses.  In~\cite{\vogtsurvey} these complexes are denoted $K_W$, and we will follow that notation here. To prove that the entire complex $K_n$ is contractible, Culler and Vogtmann prove that, for any finite set $W,$ $K_W \simeq K_n$.  They then find a set of conjugacy classes such that $K_W$ is the star of a single marked rose and therefore contractible.  Putting these two facts together we have,

\begin{lemma}[Culler-Vogtmann~\cite{\cv}]\label{L:kw.contractible}
For any finite set $W \subseteq \C$, $K_W$ is contractible.
\end{lemma}

Proposition~\ref{P:intersections} will be proved by finding a deformation retraction from $\bigcap \RR_{\Sigma_i}$ to a suitable $K_W$.  We begin by studying of the behavior of the norm with respect to Whitehead automorphisms.  For us, the traditional Whitehead automorphisms are less convenient to work with than a slightly modified version, given in~\cite{\hoare}.  This is because the effect of an automorphism on the star graph of a set of conjugacy classes (defined below), is easier to describe using this definition rather than the classical definition of Whitehead automorphism.

\begin{definition}\label{D:whitehead}  For basis $X$, and subset $A \subseteq X \cup X^{-1}$ for which there is a letter $a \in X \cup X^{-1}$ such that $a\in A$ but $a^{-1} \notin A$, the automorphism mapping $a$ to $a^{-1}$ and whose action on $X \cup X^{-1} -\{a, a^{-1}\}$ is given by, 
\begin{equation}
    \begin{cases}
        x \mapsto axa^{-1}, &\text{if } x \in A \text{ and } x^{-1} \in A;\\
        x \mapsto xa^{-1}, &\text{if } x \in A \text{ and } x^{-1} \notin A;\\
        x \mapsto ax, &\text{if } x \notin A \text{ and } x^{-1} \in A;\\
        x \mapsto x, &\text{if } x \notin A \text{ and } x^{-1} \notin A\\
    \end{cases}
\end{equation}
will be called a {\it Whitehead automorphism} and will be denoted by $(A,a)$.  
\end{definition}

\begin{warning}  This definition differs from the classical Whitehead automorphism by the fact that classical Whitehead automorphisms fix $a$.  This is the only difference, but it allows us to prove,
\end{warning}

\begin{lemma}
The set of Whitehead automorphisms given in Definition~\ref{D:whitehead} generate the group $\aut$.
\end{lemma}
\begin{proof}
The Neilson automorphisms generate $\aut$,~\cite[Theorem 3.2]{\mks}.  It is straightforward to write any Neilson automorphism as a product of the Whitehead automorphisms of Definition~\ref{D:whitehead}.
\end{proof}

Another important fact about Whitehead automorphisms of this type is a corollary to the Peak Reduction lemma, proved by Hoare~\cite[Lemma 3]{\hoare} for the definition of Whitehead automorphism used here.  To simplify notation, for $\psi \in \aut$, $x\in \fn$ and $W = \{w_1, \ldots w_k\} \subseteq \C$, we omit the parenthesis and write $\psi x$ for $\psi(s)$ and $psi W$ for $\psi(W)$.

\begin{lemma*}[Peak Reduction: Hoare~\cite{\hoare}]\label{L:whitehead}
Fix a basis $X$ of $F_n$ and finite set $W \subseteq \C$.  If there is an automorphism $\psi \in \aut$ such that $\|X\|_W \geq \|X\|_{\psi W}$ then $\psi$ can be written as a product of Whitehead automorphisms, $\psi = \tau_1\tau_2\cdots\tau_k$ such that,
\begin{multline}\label{E:whitehead}
    \|X\|_W > \|X\|_{\tau_k W} > \|X\|_{\tau_{k-1}\tau_k W} > \cdots >\|X\|_{\tau_l\tau_{l+1}\cdots\tau_k W} = \\
    \|X\|_{\tau_{l-1}\tau_l\tau_{l+1} \cdots \tau_k W} = \cdots = \|X\|_{W \cdot\psi}.
\end{multline}
\end{lemma*}

\begin{corollary}
If there is any automorphism $\psi$ such that $\|X\|_W > \|X\|_{\psi W}$, then there is a Whitehead automorphism reducing the norm.
\end{corollary}

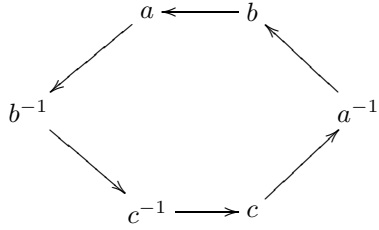
\begin{figure}
    \[
\xymatrix{
 & a \ar[dl]& b\ar[l] &\\
b^{-1} \ar[dr] & & & a^{-1} \ar[ul] \\
 & c^{-1} \ar[r] & c \ar[ur] &
}
\]
        \caption{$S_W(X)$ for $X = \{a,\;b,\;c\}, W=\{aba^{-1}b^{-1}c,\; c^{-1}\}$}\label{F:star.graph}
\end{figure}

The star graph of $W$ with respect to $X$ will allow us to study the behavior of $\|\cdot\|_W$ with respect to Whitehead automorphisms.  Recall that star graph of $W \subseteq \C$ with respect to the basis $X$ is the graph with vertex set $X \cup X^{-1}$ and with a directed edge from $x$ to $y^{-1}$ for every time the subword $xy$ appears among the conjugacy classes in $W,$ viewed as cyclic words in the alphabet $X \cup X^{-1}$;  see Figure~\ref{F:star.graph}.  The star graph of $W$ with respect to $X$ will be denoted by $S_W(X),$ or by $S_W(\rho)$ if we are thinking of $X$ as a set of labels on the marked rose $\rho.$

To prove the Peak Reduction Lemma, Hoare describes a 3-step process for constructing $S_{\tau W}(X)$ from $S_W(X)$ for a Whitehead automorphism $\tau$. If $\tau = (A,a)$ then Hoare's 3 steps are:

\begin{enumerate}
    \item Add two new vertices, $\alpha, \bar\alpha$.  Replace every edge going from a vertex in $A$ to a vertex in $A'$ (the complement of $A$) by a pair of edges, one from the vertex in $A$ to $\alpha$ and another from $\bar\alpha$ to the vertex in $A'$.  Replace every edge going from a vertex in $A'$ to a vertex in $A$ by a pair of edges, one from the vertex in $A'$ to $\bar\alpha$ and another from $\alpha$ to the vertex in $A$. 

    \item  Switch the letter $a$ with $\alpha$, and $a^{-1}$ with $\bar\alpha.$

    \item  Do the reverse of (1), reconnecting edges incident to $\alpha$ and $\bar\alpha$ according to the cyclic words in $W$ that produced them.
\end{enumerate}

Because $\|X\|_W$ is the number of edges in $S_W(X)$, this process gives the following procedure for calculating the effect of a Whitehead automorphism on the norm.  Consider the Whitehead automorphism $\tau = (A,a)$.  Draw a circle $C$ in the plane and immerse $S_W(X)$ in the plane such that each vertex of $A$ lies inside the circle, each vertex of $A'$ lies outside the circle, and such that $\#(S_W(X) \cap C)$ is minimal over all such immersions.  Then,
\begin{equation}\label{E:norm.diff}
    \|X\|_W - \|X\|_{\tau W}= val(a) - \#(S_W(X) \cap C),
\end{equation}
where $val(a)$ is the valence of the vertex of $S_W(X)$ corresponding to $a$.  

When $W = W_\Sigma$ is the set of boundary classes of a marked surface, we will be more concerned with $\|X\|_W - \|\tau^{-1}X\|_W$, because if $X$ is the set of labels on the edges of a marked rose $\rho$, then $\tau^{-1}X$ is the set of labels on $\rho \cdot \tau$.  Equation~\ref{E:norm.diff} will suffice because, 
\begin{equation}\label{E:norm}
    \|X\cdot\tau\|_W = \|X\|_{W\cdot\tau}.
\end{equation}
Equation~\ref{E:norm} follows from the observation that if $\tau(a) = x_1x_2\cdots x_k$ is an expression for $\tau(a)$ in terms of the basis $X$, $a = \tau^{-1}(x_1)\tau^{-1}(x_2)\cdots\tau^{-1}(x_k)$ is an expression for $a$ in terms of the basis $\tau^{-1}X$.  The interpretation of this observation in terms of star graphs is,
\begin{equation}\label{E:st.graph}
    S_W(\tau^{-1}X) \approx S_{\tau W}(X).
\end{equation}

\begin{lemma}\label{L:stargraph.cycles}  Let $W$ be a finite set of conjugacy classes of $\fn$.  Suppose that for some basis $X$, $S_W(X)$ is a cycle. Then $\|X\|_W$ is minimal over all bases of $F_n$, and if $Y$ is another basis with $\|Y\|_W = \|X\|_W$, then $S_W(Y)$ is also a cycle.  
\end{lemma}

\begin{proof} Since $S_W(X)$ is a cycle, any circle separating some generator from its inverse must intersect at least two edges of the graph.  Since all vertices have valence $2$, equations~\ref{E:norm.diff} and~\ref{E:st.graph} imply that no Whitehead automorphism can take $X$ to a basis that reduces the sum of the lengths of the minimal representatives for the classes in $W$.  But, if there is any automorphism reducing the sum of the lengths of the conjugacy classes in $W$, the Peak Reduction Lemma implies that a Whitehead automorphism reduces the length.  Thus, $\|X\|_W$ must be minimal over all bases for $F_n$.  

Now suppose that $Y$ is another basis with $\|Y\|_W = \|X\|_W$.  Since $\aut$ acts transitively on bases of $\fn$, we may choose $\psi \in \aut$ with $Y = \psi^{-1} X$.  By the Peak Reduction Lemma and Equation~\ref{E:st.graph}, there is a sequence of Whitehead automorphisms, $\tau_1,\ldots,\tau_l$ with $\psi = \tau_1,\ldots,\tau_l$ and $\|\tau_i^{-1}\tau_{i+1}^{-1} \cdots \tau_l^{-1} X\|_{W} = \|X\|_W$ for $i = 1 \ldots l$.  Thus, without loss of generality, we may assume that $\psi = \tau$ is a Whitehead automorphism.

We now use Hoare's method to construct the star graph $S_{\tau W}(X) \approx S_W(\tau^{-1}X)$.  Since $\|\tau^{-1}X\|_{W}=\|X\|_W$, the the circle separating $A$ from $A'$ in the star graph must intersect only two edges of the graph.  Otherwise the norm would increase.  In this case, the subgraph spanned by the vertices of $A$ is a simple path, and the same is true for the subgraph spanned by $A'$.  Step (1) of Hoare's procedure produces a graph consisting of two disjoint cycles, one containing $\alpha$ and the other containing $\bar\alpha$.  Step (2) keeps $\alpha$ and $\bar\alpha$ in separate cycles, but they may switch cycles.  Step (3) breaks these two cycles at $\alpha$ and $\bar\alpha$, and reconnects the ends of the resulting line segments to form a cycle, which is $S_{\tau W}(X)$.  Since $S_W(Y) = S_W(\tau^{-1}X) \approx S_{\tau W}(X)$, $S_W(Y)$ is a cycle.
\end{proof}

We will use Lemma~\ref{L:stargraph.cycles} to analyze which marked graphs can be drawn in a particular surface $\ms$.  To begin with, note that if $\Sigma$ has genus $g$ and $s$ punctures, and if
\begin{eqnarray*}
    X = \{a_1,b_1,a_2,b_2,\ldots,
    a_{2g},b_{2g},c_1,\ldots,c_{s-1}\}
\end{eqnarray*}
is a standard, geometric basis of $F_n = \pi_1(\Sigma)$ then,
\begin{eqnarray*}
    W_{\ms} & = & \{[a_1,b_1][a_2,b_2]\cdots[a_{2g},b_{2g}]c_1
    \cdots c_{s-1},\,c_1^{-1},\,c_2^{-1},\,\ldots,\,c_{s-1}^{-1}\}.
\end{eqnarray*} 
Thus, $S_W(X)$ is a cycle, and as a corollary to Lemma~\ref{L:stargraph.cycles} we have,

\begin{corollary}\label{C:min.rose}  If $W$ is the set of boundary classes of a surface $\Sigma$, then $$\min_\rho\|\rho\|_W = 2n$$ and $St_W(\rho')$ is a cycle for any rose $\rho'$ minimizing $\|\cdot\|_W$.
\end{corollary}

The next two lemmas characterize the marked graphs that lie in $\RR_{\ums}$.

\begin{lemma}\label{L:graph.in.surface}  The marked graph $\Gamma = (\Gamma,g)$ can be drawn in $(\Sigma,s)$ if and only if the set of reduced edge cycles of $\Gamma$ representing $W_\Sigma$ traverses each edge of $\Gamma$ exactly once in each direction.
\end{lemma}

\begin{proof}

Suppose that $\mg$ can be drawn in $\ms.$  Cutting $\Sigma$ along $\Gamma$ produces a collection of punctured disks, one for each puncture.  The oriented boundaries of these disks correspond the the conjugacy classes of the boundary of $\Sigma$.  Together they traverse each edge of $\Gamma$ once in each direction.  Thus if $\mg$ can be drawn in $\ms$ then the set of reduced edge cycles of $\Gamma$ representing the boundary classes of $\ms$ traverses each edge exactly once in each direction.

For the converse, we will first construct another marked surface $\Sigma^\prime$ whose set of oriented boundary classes is also $W_\Sigma$ and in which $\Gamma$ can be drawn.  We then use this surface to draw $\Gamma$ in $\Sigma$.  Suppose that the set of reduced edge cycles in $\Gamma$ that represents $W_\Sigma$ traverses each edge exactly once in each direction.  These edge cycles give $\Gamma$ a unique ribbon structure $\OO$ such that the set of boundary classes of $\mrg$ is exactly $W_\Sigma$.  

The surface $\Sigma' = |\mrg|$ is is homeomorphic to $\Sigma$ because they are both orientable surfaces with the same number of punctures and the same fundamental group.  The set of boundary classes of $\Sigma'$ is $W_\Sigma$.  Let $f:\Sigma^\prime \rightarrow\Sigma$ be any homeomorphism that preserves the labels of the punctures.  Now, $f\circ i$ embeds $\Gamma$ into $\Sigma$ as a strong deformation retract, but this embedding may not induce the same marking as $g$.  That is to say, the following diagram may not commute up to homotopy,

\[
\xymatrix@C-15pt@R-15pt{
\Gamma\ar@{^{(}->}[rr]^i&  & \Sigma'\ar[dd]^f\\
& & \\
R_0\ar[uu]^g \ar[rr]_s & & \Sigma.\\
}
\]

However, the outer automorphism given by $f_*\circ i_*\circ g_* \circ (s_*)^{-1}$ stabilizes $W_\Sigma$.  By a theorem of Zieschang~\cite[Theorem 5.15.3]{\zies}, it is induced by an element $\theta$ in the orientatin preserving mapping class group of $\Sigma$.  Embedding $\Gamma$ into $\Sigma$ by $\theta^{-1}\circ f\circ i$ gives the same marking as $g$.  Thus, $\Gamma$ can be drawn in $\Sigma$.
\end{proof}

\begin{lemma}\label{L:graph.in.kmin}  Let $\ms$ be a homotopy marked surface and $\rho$ a marked rose.  Then $\rho\in K_{W_\Sigma}$ if and only if $\rho \in \RR_\Sigma$.
\end{lemma}

\begin{proof} By Corollary~\ref{C:min.rose}, the minimal value of $\|\cdot\|_W$ is $2n$.  If $\rho \in \RR_\Sigma$ then $\rho$ can be drawn in $\Sigma$.  By cutting $\Sigma$ along $\rho$, we see that $\|\rho\|_{W_\Sigma} =2n$ so that $\rho\in K_{W_\Sigma}.$

Conversely suppose that $\rho \in K_{W_\Sigma}$.  Since $\rho$ minimizes $\|\cdot\|_{W_\Sigma}$, Corollary~\ref{C:min.rose} implies that the star graph $S_{W_\Sigma}(\rho)$ must be a cycle.  Therefore, each label in $\rho$ appears exactly once with exponent $+1$ and once with exponent $-1$ in the minimal expressions for conjugacy classes of $W_\Sigma$ in terms of a set of labels of $\rho$.  This means that the set of reduced edge cycles in $\rho$ that represents $W_\Sigma$ traverses each edge of $\rho$ exactly twice, once in each direction.  By Lemma~\ref{L:graph.in.surface}, $\rho$ can be drawn in $\ms$. 
\end{proof}

This lemma says that the roses in $\RR_\Sigma$ coincide with the roses in $K_{W_\Sigma}$.  Since $K_{W_\Sigma}$ is the union of the stars of its roses, $\RR_\Sigma \subseteq K_{W_\Sigma}.$  To find a graph in $\RR_\Sigma$ lying near a particular graph in $K_{W_\Sigma} - \RR_\Sigma$, we use the following lemma.

\begin{lemma}\label{L:min.forest}   If $\Gamma = \mg \in K_{W_\Sigma}$ then there exists a (possibly empty) forest $\Phi_\Sigma(\Gamma)$ such that for any other forest $\Phi \subseteq \Gamma$,

$$ \Gamma / \Phi \in \RR_\Sigma \Longleftrightarrow \Phi\supseteq\Phi_\Sigma(\Gamma).$$
\end{lemma}

\begin{proof}  Let $\Phi_\Sigma(\Gamma)$ be the subgraph of $\Gamma$ consisting of all the edges of $\Gamma$ that are not traversed exactly once in each direction by the set of reduced edge paths representing the boundary classes of $\Sigma$.  Since $\Gamma \in K_{W_\Sigma}$, and $K_{W_\Sigma}$ is the union of the stars of its roses, there is a maximal tree $T$ in $\Gamma$ such that $\Gamma / T$ is a rose in $K_{W_\Sigma}$.  By Lemma~\ref{L:graph.in.kmin} this rose is in $\RR_\Sigma$, so it can be drawn in $\Sigma$.  Therefore, every edge of $\Gamma - T$ is traversed exactly once in each direction by the set of conjugacy classes in $W_\Sigma$.  This means that $\Phi_\Sigma(\Gamma)\subseteq T$, so that  $\Phi_\Sigma(\Gamma)$ is a forest.

Now, given any forest $\Phi$ in $\Gamma$, Lemma~\ref{L:graph.in.surface} implies that  $\Gamma / \Phi$ can be drawn in $\Sigma$ exactly when the boundary cycles traverse each edge of $\Gamma / \Phi$ once in each direction.  This happens exactly when $\Phi_\Sigma(\Gamma) \subseteq \Phi$.
\end{proof}

These lemmas would allow us, at this time, to define a retraction from $K_{W_\Sigma}$ to $\RR_\Sigma$ by taking a graph $\Gamma \in K_W$ to $\Gamma / \Phi_\Sigma(\Gamma)$ thus proving the contractibility of $\RR_\Sigma$ for a marked surface $\Sigma$.  But, we will postpone this until it is covered by the proof of contractibility for arbitrary simplices.  For higher dimensional simplices, we need a set of conjugacy classes that captures the properties of a graph that can be drawn in several different surfaces.  This set emphasizes a conjugacy class according to the number of the surfaces in question of which it is a boundary.  We start by describing some general properties of collections of finite sets of conjugacy classes of $\fn$.  For the proof of Proposition~\ref{P:intersections}, we will specialize to the case that the sets of conjugacy classes are actually the boundary classes of marked surfaces.

\begin{definition}\label{D:boundary.classes}  For a collection $\sigma = \{W_0, \ldots, W_k \}$ of finite sets of conjugacy classes of $\fn$, define
$$W_\sigma := \{[\alpha_1]^{n_1},\ldots,[\alpha_l]^{n_l}\},$$
 where $\bigcup_{i=0}^k W_{i} = \{[\alpha_1],\ldots,[\alpha_l]\}$, and $n_j$ is the number of times that the conjugacy class $[\alpha_j]$ appears in the $W_{i}$.  
\end{definition}

Note that $[\alpha]$ and $[\alpha^{-1}]$ both may appear in $W_\sigma$.  We use the letter $\sigma$ for the set $\{W_0, \ldots, W_k \}$ because this definition will be applied to a simplex $\sigma = \langle \Sigma_0, \ldots, \Sigma_k \rangle$ of $\NN$, with $W_i = W_{\Sigma_i}$. We will use the notation $W_\sigma$ in this situation as well.

\begin{lemma}\label{L:norm.sum}  
For $\sigma$ and $W_\sigma$ as above, let
\begin{eqnarray*}
    A & = &\min_{\rho}\|\rho\|_{W_\sigma}, \\
    A_i & = & \min_\rho\|\rho\|_{W_i}.
\end{eqnarray*}
Then $A = A_0 + \cdots +A_k$ if and only if $\bigcap _{i=0}^k K_{W_i} \neq \emptyset$.
\end{lemma}

\begin{proof}  
Suppose that $W_\sigma = \{w_1^{n_1},\ldots,w_l^{n_l}\}$.  For any marked rose $\rho$
\begin{equation}\label{E:norm.expressions}
\|\rho\|_{W_\sigma} = \sum_{i=0}^l n_i\|\rho\|_{\{w_i\}} = \sum_{j=0}^k \|\rho\|_{W_j}.
\end{equation}

Choose any marked rose $\rho_1$ with $\|\rho_1\|_{W_\sigma} = A$.  Now, $\rho_1$ may not minimize every $\|\cdot\|_{W_i}$, so
\begin{equation}\label{E:sum.less} 
    A_0 + \cdots +A_k \leq \sum_{j=0}^k \|\rho_1\|_{W_j} = \|\rho_1\|_{W_\sigma} = A,
\end{equation}
where the first equality comes from equation~\ref{E:norm.expressions}.

Now, if $\bigcap _{i=0}^k K_{W_i} \neq \emptyset$ then there is a single marked rose, $\rho_2$ with $\|\rho_2\|_{W_i} = A_i$ for all $i$.  Thus,
\begin{equation*}
    A\leq \|\rho_2\|_{W_\sigma} = \sum_{i=0}^k\|\rho_2\|_{W_i} = A_0 +\cdots+A_k.
\end{equation*}
Together with~\ref{E:sum.less} this implies that $A= A_0 +\cdots+A_k$.

Conversely if $A = A_0 +\cdots+A_k$, then using the $\rho_1$ from above we have, 
\begin{equation}\label{Eq:a.sum}
    A_0 +\cdots+A_k = A = \|\rho_1\|_{W_\sigma} = \sum_{j=0}^k \|\rho_1\|_{W_j}.
\end{equation}
Again the last equality comes from equation~\ref{E:norm.expressions}.  Now,  $\|\rho_1\|_{W_i}\geq A_i$, so by~\ref{Eq:a.sum} $\|\rho_1\|_{W_i} = A_i$ for each $i$.  Hence, $\rho_1\in K_{W_i}$ for each i, and $\bigcap K_{W_i} \neq \emptyset$.
\end{proof}

Changing the viewpoint slightly gives the following corollary.

\begin{corollary}\label{C:nonempty.kmin}  For any finite collection of finite sets of conjugacy classes, $\sigma = \{W_0, \ldots, W_k \}$, $K_{W_\sigma} = \bigcap_{i=0}^{k} K_{W_i}$ if the right hand side is nonempty.
\end{corollary}

The final lemma we need for the proof of Proposition~\ref{P:intersections} is the Poset Lemma of~\cite{\quillen}.

\begin{lemma*}[Poset Lemma]\label{L:poset}
Let $f:P \to P$ be a poset map from the poset $P$ to itself such that $p \leq f(p)$ for all $p \in P$.  Then $f$ induces a deformation retraction from the geometric realization of $P$ to the geometric realization to its image, $f(P)$.
\end{lemma*}


We are now in a position to prove Proposition~\ref{P:intersections}.

\begin{proof}[Proof of Proposition~\ref{P:intersections}]  
Let $\sigma = \langle \Sigma_0, \ldots, \Sigma_k \rangle$ be a simplex of $\NN_n$.  Denote by $\RR_\sigma$ the intersection,
$$\RR_\sigma := \bigcap_{i=0}^k \RR_{\ums[i]},$$
and let $W_\sigma$ be the set of conjugacy classes given by Definition~\ref{D:boundary.classes}. Since $ \langle \Sigma_0 \ldots \Sigma_k \rangle$ is a simplex of $\NN_n$, $\RR_\sigma$ contains a rose.  To simplify the notation, set $W_i := W_{\Sigma_i}$.  Since $\RR_{\Sigma_i} \subseteq K_{W_i}$, $\RR_\sigma \subseteq \bigcap K_{W_i}.$  Hence, $\bigcap K_{W_i} \neq \emptyset$ and by Corollary~\ref{C:nonempty.kmin}, $\bigcap K_{W_i} = K_{W_\sigma}$.  We will define a deformation retraction $K_{W_\sigma} \rightarrow \RR_\sigma$ by collapsing in each graph the minimal forests that takes that graph to a graph in $\RR_\sigma$.

Let $\Gamma \in K_{W_\sigma}$.  Since $\bigcap K_{W_i} = K_{W_\sigma}$, Lemma~\ref{L:min.forest} gives the minimal forest $\Phi_{\Sigma_i}$ collapsing $\Gamma$ to a graph in $\RR_{\Sigma_i}$.   Set 
$$\Phi_\sigma(\Gamma) := \Phi_{\Sigma_0}(\Gamma) \cup \cdots \cup \Phi_{\Sigma_k}(\Gamma).$$
Since $\Gamma \in K_{W_\sigma}$, there is a spanning tree $T$ collapsing $\Gamma$ to a rose, $\Gamma / T \in K_{W_\sigma}$.  By Corollary~\ref{C:nonempty.kmin}, $\Gamma / T \in K_{W_i}$ for each $i.$  So, by Lemma~\ref{L:graph.in.kmin}, $\Gamma / T \in \RR_{\Sigma_i}$ for each $i$ and therefore, $\Phi_{\Sigma_i}(\Gamma) \subseteq T$ for each $i$.  Therefore, $\Phi_\sigma(\Gamma) \subseteq T$.  Since $T$ is a tree, $\Phi_\sigma(\Gamma)$ is a forest.

Now we define a map $r$ from the vertex set of $K_{W_\sigma}$ to the vertex set of $\RR_\sigma$ by $r(\Gamma) = \Gamma / \Phi_\sigma(\Gamma)$.  We claim that $r$ induces a simplicial map,
\begin{eqnarray*}
 r:K_{W_\sigma} &\longrightarrow& \RR_\sigma.
\end{eqnarray*}
It will suffice to show that $r$ takes adjacent vertices to the same vertex or adjacent vertices because both $K_{W_\sigma}$ and $\RR_\sigma$ are determined by their 1-skeletons.  To to this, suppose that $\Gamma_1$ and $\Gamma_2$ represent adjacent vertices in $K_{W_\sigma}$.  By possibly switching the names of the graphs, we can write 
$\Gamma_2 = \Gamma_1/\Phi$
for some forest $\Phi$.  If $\Phi \subseteq \Phi_\sigma(\Gamma_1)$, then $r(\Gamma_1)=r(\Gamma_2)$.  If $\Phi \not\subseteq \Phi_\sigma(\Gamma)$ then $r(\Gamma_1) \neq r(\Gamma_2)$.  The diagram below represents a small portion of $K_n$ in this case, with edges represented by arrows.

\[
\xymatrix@C-5pt@R-5pt{
r(\Gamma_2)& &  r(\Gamma_1)\ar[ll]_{collapse \;\Phi''}\\
& & \\
\Gamma_2\ar[uu]^{collapse \;\Phi_\sigma(\Gamma_2)}&  & \Gamma_1
    \ar[ll]_{collapse \;\Phi}
    \ar[uull]|-{collapse \;\Phi^\prime}
    \ar[uu]_{collapse \;\Phi_\sigma(\Gamma_1)}
}
\]

In order to show that $r$ takes $\Gamma_1$ and $\Gamma_2$ to adjacent vertices, as the diagram suggests, we need justify that there is such a forest $\Phi''$, as indicated in the diagram.  The forest $\Phi''$ is constructed as follows.  Let $ \Phi^\prime = \Phi_\sigma(\Gamma_1)\cup\Phi.$ 
Then $\Phi'$ is the subforest of $\Gamma_1$ such that 
$r(\Gamma_2) = \Gamma_1/\Phi^\prime$.
If $\Phi^{\prime\prime}$ is the subgraph of $r(\Gamma_1)$ consisting of the images of the edges in 
$\Phi^\prime - \Phi_\sigma(\Gamma_1)$
then $\Phi^{\prime\prime}$ is a forest and 
$r(\Gamma_1)/\Phi^{\prime\prime} = r(\Gamma_2)$.
Therefore, $r(\Gamma_1)$ and $r(\Gamma_2)$ are adjacent, proving that $r$ induces a simplicial map.

That $r$ is a retraction follows from the fact that, if $\Gamma\in \RR_\sigma$ then $\Phi_\sigma(\Gamma) =\emptyset$.  Also, by Lemma~\ref{L:graph.in.surface}, the image of $r$ is contained in $\RR_{\Sigma_i}$ for each $i$.  Therefore, $r(K_{W_\sigma}) \subseteq \RR_\sigma$.  To see that $r$ is a deformation retraction, we will use the Poset lemma.  Partially order the vertices of $K_{W_\sigma}$ by, $\Gamma_1 < \Gamma_2$ if $\Gamma_1$ can be collapsed to $\Gamma_2.$  Then, $K_{W_\sigma}$ is the geometric realization of the poset of its vertices under this partial order.  With respect to this partial order, $r$ has the property that $\Gamma \leq r(\Gamma)$.  The Poset lemma implies that $r$ is a deformation retraction.  Since $K_{W_\sigma}$ is contractible, this finishes the proof.
\end{proof}

By Lemma~\ref{L:nerve.contract}, Proposition~\ref{P:intersections} proves that $\NN_n \simeq K_n$.  Since $K_n$ is contractible~\cite{\cv}, so is $\NN_n$.  We record this as,

\begin{theorem}\label{T:contractible}
$\NN_n$ is contractible for all $n$.
\end{theorem}

\section{Simplex stabilizers}
As mentioned before, the action of $\out$ on $\kn$ gives an action of $\out$ on $\NN_n$.  We prove that the vertex stabilizers for this action are extended mapping class groups.  In this section we deal strictly with punctured surfaces without boundary and their extended mapping class groups, in which mapping classes are allowed to reverse orientation.  As usual, let $\ms$ be a homotopy marked surface with corresponding ribbon graph subcomplex $\RR_{\ums}$.  The identification of $\pi_1(\Sigma)$ with $\fn$ given by the marking $s$ induces a homomorphism from the extended mapping class group of $\Sigma$ to $\out$.  This homomorphism is defined by sending a homeomorphism of $\Sigma$ to the outer automorphism of $\pi_1(\Sigma)$ that it represents.  By Zieschang's theorem~\cite[Theorem 5.15.3]{\zies}, this homomorphism is injective, and its image is the subgroup of $\out$ consisting of outer automorphisms that take $W_\Sigma$ to $W_\Sigma$ or $(W_\Sigma)^{-1}$.   Denote this subgroup by $\mcg\ms$ and denote the image of the orientation preserving subgroup by $MCG\ms$.  Note that these subgroups depend on the marking $s$.   Let $\Stab(\RR_{\ums})$ be the subgroup of ${\rm Out}(F_n)$ stabilizing $\RR_{\ums}$ setwise, so $\Stab(\RR_{\ums})$ is the stabilizer of the vertex of $\NN_n$ that corresponds to $\ums$.

\begin{theorem}\label{T:vertex.stabilizers}  $\Stab(\RR_{\ums})= \mcg\ms.$
\end{theorem}

\begin{proof}

The fact that $\Stab(\RR_{\ums}) \supseteq \mcg\ms$ follows from the fact that if $\psi \in \mcg\ms$, then $\ms \cdot \psi = (\Sigma,s \circ |\psi|) \sim_\pm \ms.$  Therefore, $\RR_{\ums} \cdot \psi = \RR_{\ums}$ showing that $\Stab(\RR_{\ums}) \supseteq \mcg\ms$.

On the other hand, suppose that $\psi \in \Stab(\RR_{\ums})$.  Then, 
$$\RR_{\ums} = (\RR_{\ums})\cdot\psi = \RR_{[\Sigma,s\circ|\psi|]}.$$
By Proposition~\ref{P:equivalence}, $\ums =[\Sigma,s\circ|\psi|]$, so there is a homeomorphism $h : \Sigma \to \Sigma$ that makes the following diagram commute up to homotopy,

\[
\xymatrix{
R_0 \ar[r]^s \ar[d]_{|\psi|} & \Sigma \ar[d]^h\\
R_0 \ar[r]_s & \Sigma.
}
\]

Now, $h$ takes the boundary classes $W_\Sigma$ to $W_\Sigma$ or $(W_\Sigma)^{-1}$.  Thus, $\psi$ does also.  By Zieschang's theorem, $\psi\in\mcg(\Sigma)$.
\end{proof}

To describe the stabilizer of a higher dimensional simplex, we use a certain kind of stabilizer of a set of conjugacy classes of $\fn$, as studied by McCool in~\cite{\mccooltwo}.  Following the definitions there, we consider ordered $m$-tuples of conjugacy classes in $\fn$,
$$(w_1,\ldots,w_m).$$
The symmetric group, $S_m$ acts on the set of $m$-tuples by permuting the coordinates.  The {\em inverting operations}, $\tau_1,\ldots,\tau_m$ act on the set of $m$-tuples by,
\begin{equation*}
\tau_i(w_1,\ldots,w_i,\ldots,w_m) = (w_1,\ldots,w_i^{-1},\ldots,w_m)
\end{equation*}
The group $S_m$ together with the $\tau_i$'s generate the subgroup $\Omega_m \cong S_m\wr \Z_2$, of permutations of the set of $m$-tuples of conjugacy classes of $\fn$ known as the extended symmetric group.  The group $\out$ also acts on the set of $m$-tuples of conjugacy classes by acting individually on the coordinates.  In~\cite{\mccooltwo}, McCool makes the following definition in the setting of $\aut$, but which we will use in the setting of $\out$.

\begin{definition}
For $U$ an $m$-tuple of conjugacy classes and subgroup $G \leq \Omega_m$, define the subgroup $\A_{U,G}$ of $\out$ by
$$\A_{U,G} := \{\theta \in \out | \theta U \in GU \},$$
where $GU = \{gU | g \in G \}$.
\end{definition}

For a simplex $\sigma$ of $\NN_n$, let $G_\sigma$ denote the stabilizer of the simplex $\sigma$ of $\NN_n$.  If $\sigma = v$ is the vertex corresponding to the marked surface $\ms$, then $G_v = \mcg\ms = \A_{U,G}$ where $U$ is the $m$-tuple of boundary classes (in any order) of a marked surface $\ms$, and $G \leq \Omega_m$ is the subgroup generated by $S_m$ together with the extended permutation $\tau_1\tau_2\cdots\tau_m$.  To describe $G_\sigma$ for a higher dimensional simplex, we make the following,

\begin{definition}\label{D:simplex.stab}
Let  $\sigma = \langle \Sigma_0,\ldots,\Sigma_k \rangle$ be a simplex of $\NN_n$, and let $U_i$ be the $m_i$-tuple of boundary classes of $\Sigma_i$ (again in any order).  Denote by $U_\sigma = (U_0,\ldots,U_k)$ the $m = (m_0 + \cdots + m_k)$-tuple constructed by deleting the appropriate pairs of parentheses.  Define $H_\sigma \leq \Omega_m$ be the subgroup generated by extended permutations of the following types,
\begin{enumerate}
    \item $\alpha \in S_m$ such that there is a permutation $\lambda \in S_k$ such that, for each $i$, $\alpha$ takes $U_i$ to $U_{\lambda(i)}$, possibly with the entries of $U_{\lambda(i)}$ permuted.

    \item $\tau \in \Omega_m$ such that, for each $i$ $\tau U_i = U_i$ or $\tau U_i = U_i^{-1}$.
\end{enumerate}
\end{definition}

\begin{proposition}\label{P:simplex.stab} For any simplex $\sigma$ of $\NN$, $G_\sigma$ has the form $\A_{U_\sigma,H_\sigma}$.
\end{proposition}

\begin{proof}
Formally, $(1)$ can be written, $\alpha U_i \in S_{m_{\lambda(i)}} U_i$.  If $\theta \in G_\sigma$ then $\theta$ permutes the equivalence classes of the marked surfaces $\Sigma_i$.  This means that $\theta$ takes $W_{\Sigma_i}$ to $W_{\Sigma_j}$ or $W_{\Sigma_j}^{-1}$ for some $j$.  Thus, $\theta U_i \in S_{m_j} U_j \text{ or } S_{m_j} U_j^{-1}$ for some $j$.  Since no two surfaces are taken to the same surface by $\theta$, this means precisely that $\theta \in \A_{U_\sigma,G_\sigma}$ as defined above.
\end{proof}

\section{Equivariant homology of the action of $\out$ on $\NN_n$}\label{S:ss}

For a cellular action of a group $G$ on a contractible cell complex $X$, the equivariant spectral sequence for the action of $G$ on $X$ is a well-known spectral sequence that converges to a grading of the homology of $G$ (see~\cite[Chapter VII.7]{\brown}).  To describe this spectral sequence, let $M$ be any $G$-module.  Consider, for each $p$-cell $\sigma$ of $X$, the $G_\sigma$-module $\Z_\sigma$. As an additive group, $\Z_\sigma$ is isomorphic to $\Z$.  The module structure of $\Z_\sigma$ is given by having $g\in G_\sigma$ act as multiplication by $+1$ or $-1$, depending on whether $g$ preserves or reverses the orientation of $\sigma$.  The module $\Z_\sigma$ is called the {\em orientation module} of $\sigma$.  Let
$$M_\sigma := \Z_\sigma \otimes_\Z M.$$   
Fix a set $\Delta_p$ of representatives for the orbits of the $p$-cells of $X$  under the action of $G$.  The equivariant spectral sequence for the action of $G$ on $X$ takes the form:

\begin{equation}\label{E:e1}
E^1_{pq} = \bigoplus_{\sigma\in\Delta_p}H_q(G_\sigma;M_\sigma)\Rightarrow H_{p+q}(G;M).
\end{equation}

Applying this to the action of $\out$ on $\NN_n$ and any $\out$-module, $M$, gives a spectral sequence converging to $H_*(\out;M)$.  Since vertex stabilizers are extended mapping class groups and there is one orbit of vertex for each homeomorphism type of surface, the $p=0$ column of the spectral sequence consists of direct sums of the homology groups of mapping class groups.  For $p > 0$, the simplex stabilizers are given by Proposition~\ref{P:simplex.stab} and we have,
\begin{theorem}\label{T:out.ss} For any $\out$-module, $M$, there is a spectral sequence of the form
\begin{equation}\label{E:t.e1}
E^1_{pq} = \bigoplus_{\sigma\in\Delta_p}H_q(G_\sigma;M_\sigma)\Rightarrow H_{p+q}(\out;M),
\end{equation}
where $\Delta_0$ is the set of homeomorphism classes of punctured orientable surfaces with fundamental group $\fn$ and for vertex $v\in\Delta_0$ corresponding to surface $\Sigma$, the stabilizer $G_v$ is the extended mapping class group $\mcg(\Sigma)$.  Moreover, for $p>0$, each $G_\sigma$ is a generalized stabilizer of the form $\A_{U_\sigma,H_\sigma}$.
\end{theorem}

\section{Analysis of $E^\infty$}

In this section, we specialize to $\Q$ coefficients and sketch a method for using Harer's homology stability theorems for mapping class groups~\cite{\harer} to analyze the $E^\infty$ page of this spectral sequence.  For complete proofs, see~\cite{\horak}.  In Section~\ref{S:intro}, we defined the pure mapping class of a bounded surface to be isotopy classes of homeomorphisms fixing the boundary pointwise.  We keep that convention here, and extend the definitions to surfaces with boundary and punctures.  Let $\pmcg[s]{g}{r}$ denote the pure mapping class group of the surface $\Sigma$ of genus $g$ with $s$ marked points and $r$ boundary components.  Elements of $\pmcg[s]gr$ are isotopy classes of homeomorphisms of $\Sigma$ that fix the boundary and marked points pointwise.  For $r=0$ let $\fmcg[s]{g}{0}$ denote the extended mapping class group of the surface of genus $g$ with $s$ punctures, as defined in Section~\ref{S:intro}.  By restricting a mapping class in $\pmcg[s]g0$ to the surface minus the marked points, we get a map $\pmcg[s]g0 \to \fmcg[s]g0$.  If $\Sigma$ is a particular surface, we to use $P\Gamma(\Sigma)$ and $\Gamma(\Sigma)$ to denote its pure and extended mapping class groups.

Suppose that $\Sigma_0$ is a subsurface with boundary of the surface $\Sigma$, where the boundary of $\Sigma_0$ need not be contained in the boundary of $\Sigma$.  The inclusion $\Sigma_0 \into \Sigma$ induces a map $\alpha : P\Gamma(\Sigma_0) \to P\Gamma(\Sigma)$ defined by extending a homeomorphism of $\Sigma_0$ to all of $\Sigma$ by the identity.  Harer's stability theorems imply that if the genus of $\Sigma_0$ is at least $3k-1$, then $\alpha_* : H_k(P\Gamma(\Sigma_0);\Q) \to H_k(P\Gamma(\Sigma);\Q)$ is injective.

\begin{figure}
    \psfrag{a}{\huge $a$}
    \psfrag{b}{\huge $b$}
    \psfrag{c}{\huge $c$}
    \psfrag{d}{\huge $d$}
    \psfrag{e}{\huge $e$}
    \psfrag{f}{\huge $f$}
    \begin{center}
      \resizebox{2.0 in }{1.76 in}{\includegraphics{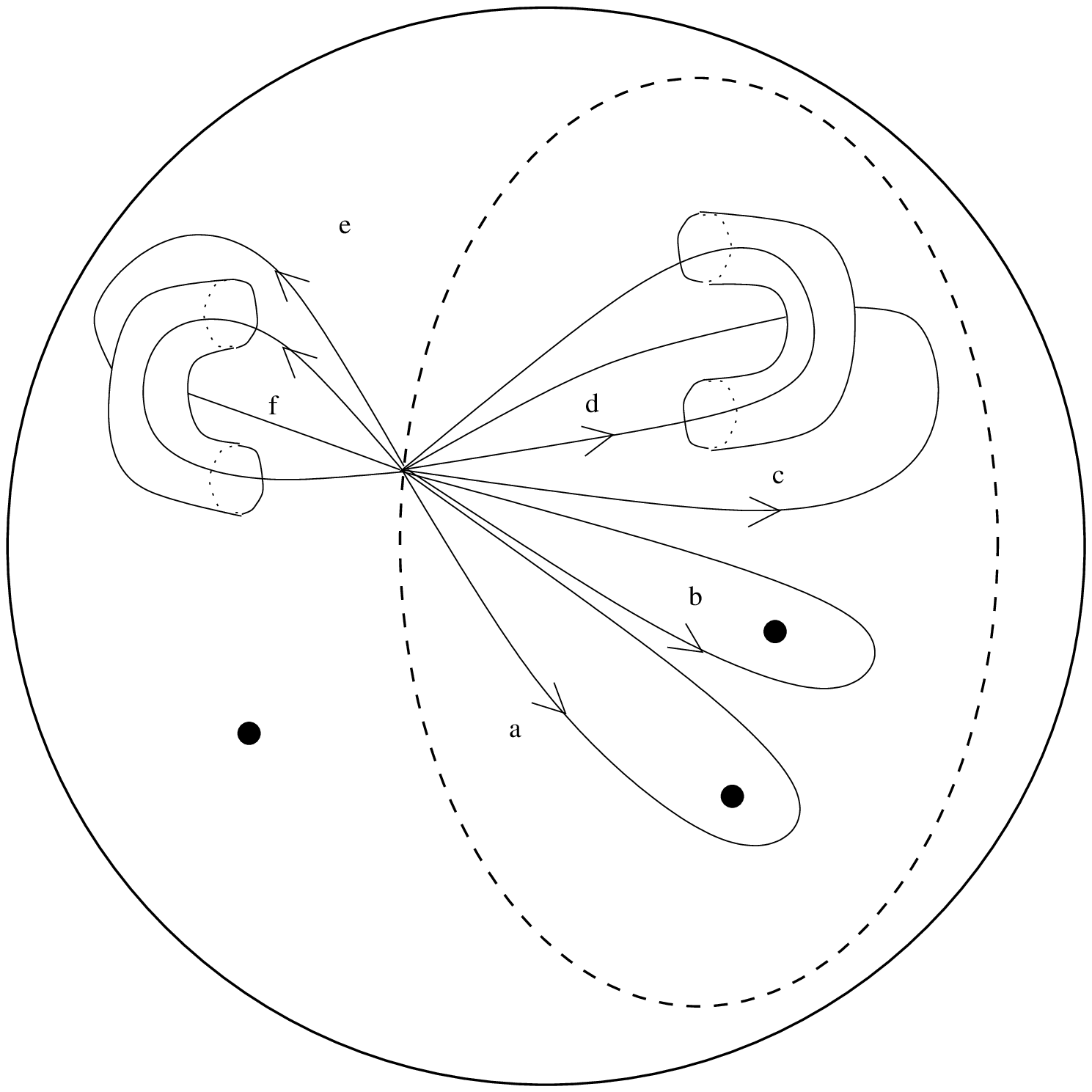}}
      \hspace{.2 in}
      \resizebox{2.0 in}{1.8 in }{\includegraphics{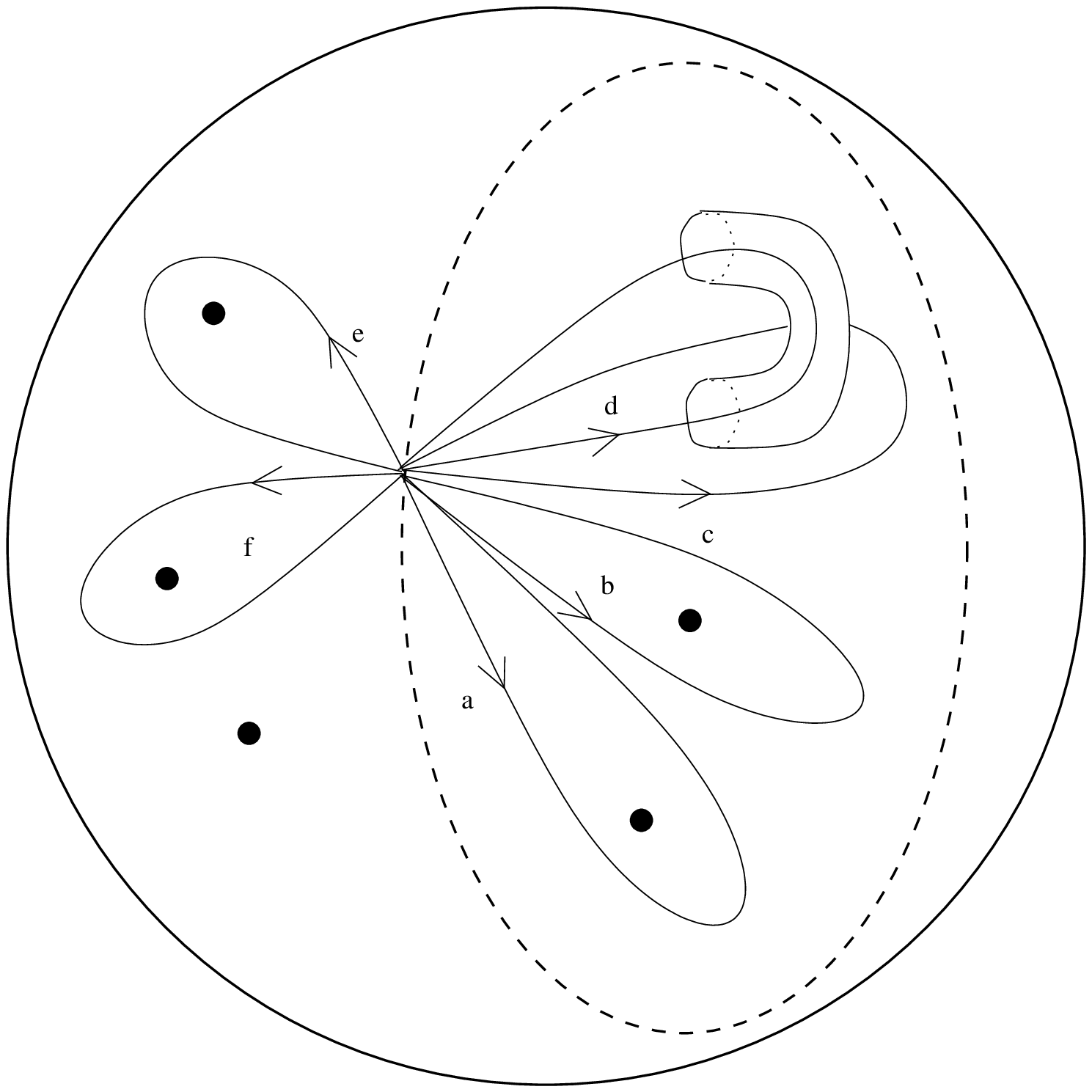}}
    \end{center}
    \caption{Markings that agree on a subsurface}\label{F:markings.agree}
\end{figure}

In order to relate these stability maps to the $d^1$ terms in the spectral sequence of Theorem~\ref{T:out.ss}, consider two marked surfaces, $\ms, (\Sigma',s')$.  Let $\rho$ and $\rho'$ be the images in $\Sigma$ and $\Sigma'$ of the marking rose.  Suppose that there are separating simple closed curves $\gamma \subset \Sigma$ and $\gamma' \subset \Sigma'$ cutting off subsurfaces $\tsigma \subset \Sigma$ and $\tsigma' \subset \Sigma'$ with the following properties, illustrated in Figure~\ref{F:markings.agree}.
\begin{enumerate}
    \item
The basepoints of the marking roses lie on $\gamma$ and $\gamma'$
    \item
No edge of the marking roses meets $\gamma$ or $\gamma'$ anywhere but at the basepoints of the roses.
    \item $\rho \cap \tsigma \simeq \tsigma$ and $\rho' \cap \tsigma' \simeq \tsigma'$.
    \item
$\tsigma$ and $\tsigma'$ are homeomorphic by a homeomorphism taking each directed edge of $\rho \cap \tsigma$ to a directed edge of $\rho' \cap \tsigma'$ with the same labelling.
\end{enumerate}

\begin{definition}\label{D:markings.agree}
If the marked surfaces $\ms$ and $(\Sigma',s')$ satisfy conditions $(1) - (4)$ above, the markings are said to {\em agree on the subsurfaces $\tsigma$ and $\tsigma'$}.
\end{definition}

In the context of marked surfaces with markings agreeing on subsurfaces, Harer's stability theorems can be used to prove,

\begin{lemma}\label{L:d1.rank}
Let $\ums$ and $[\Sigma',s']$ be two non-homeomorphic marked surfaces satisfying the following three conditions,
\begin{enumerate}
    \item  The markings $s$ and $s'$ agree on subsurfaces $\tsigma$ and $\tsigma'$ of genus $3k-1$,
    \item  $\tsigma$ contains all but one of the punctures of $\Sigma$,
    \item $\RR_{\ums} \cap \RR_{[\Sigma',s']} \neq \emptyset$.
\end{enumerate}
Let $v$ and $v'$ be the vertices of $\NN_n$ corresponding to $\Sigma$ and $\Sigma'$ and let $e$ be the edge of $\NN_n$ between $v$ and $v'$.  Then, $G_e = G_v \cap G_{v'}$ and $i_* : H_k(G_e;\Q) \to H_k(G_v;\Q)$ has rank at least
$$\dim\bigl(H_k(G_v;\Q)\bigr) - \bigl[ \dim\bigl(H_k(P\Gamma(\Sigma);\Q)\bigr) - \dim\bigl(H_k(P\Gamma(\tsigma);\Q)\bigr)\bigr].$$
\end{lemma}

By using these lemmas and Harer's stability theorems, one can find rough upper bounds on the dimensions of some parts of the $E^\infty$ page of the spectral sequence in Theorem~\ref{T:out.ss}.  For example, one can prove,

\begin{proposition}\label{P:stable.homology}
Let $k \geq 0$ and $n \geq 6k - 2$.  For $g \geq 3k-1$, let $\Sigma_g^s$ be the punctured surface of genus $g$ with $s$ punctures and with $2g+s-1=n$  (so that $\pi_1(\Sigma_g^s) \cong \fn$).  By choosing particular markings of the $\Sigma_g^s$, we may identify the vector space
$$A := \bigoplus_{g \geq 3k-1} H_k(MCG(\Sigma_g^s);\Q)$$
 with a subspace of the $E^1_{0k}$ term of spectral sequence~\ref{E:t.e1} using trivial $\Q$ coefficients. The image of $A$ in $E^\infty_{0k}$ has dimension no larger than $$\dim\bigl(H_k(\fmcg[t]{3k-1}{0};\Q)\bigr) + \sum_{\substack{g \geq 3k\\2g+s-1=n}}\bigl[\dim\bigl(H_k(\pmcg[s]{g}{0} ;\Q)\bigr) - \dim\bigl(H_k(\pmcg[s-1]{g}{0};\Q)\bigr)\bigr],$$
where $2(3k-1)+t-1=n$.
\end{proposition}

This proposition is proved by analyzing the rank of the $d^1$ map of the spectral sequence.  The strategy is to choose vertices $v_0,v_1,\cdots$, that represent the orbits of vertices that correspond to the homotopy marked surfaces $\ums[0],\ums[1],\cdots$.  Choose these representatives so that there are edges $e_g$ spanned by the vertices $v_0$ and $v_g$.  In this case, the $(e_g,v_g)$-component of the $d^1$ map is simply the map on homology induced by the inclusion $G_{e_g} \into G_{v_g}$.  The proposition is proved by choosing markings so that the markings for $v_g$ and $v_0$ agree on a subsurface of genus $3k-1$ and then using Lemma~\ref{L:d1.rank} to estimate the rank of the $(e_g,v_g)$-component of $d^1$.

\section{Acknowledgments}
This paper presents a portion of the results of my thesis, completed at Cornell University in 2003.  I would like to thank my advisor Karen Vogtmann for all of the time and support she gave me during my time at Cornell.  These ideas would not have been possible without the many useful conversations we had about this material.

\bibliographystyle{amsplain}
\bibliography{outerspace}

\end{document}